\documentclass[11pt]{article}
\usepackage{graphicx}
\usepackage{caption}
\usepackage{subcaption}
\usepackage{verbatim}

\usepackage[compress]{cite}
\usepackage{amsmath}
\usepackage{amsfonts}
\usepackage[papersize={8.5in,11in},margin=1in]{geometry}
\usepackage[bookmarks=true,         bookmarksnumbered=true,
colorlinks=true,
pdfstartview=FitV,
linkcolor=blue, citecolor=blue, urlcolor=blue,
]{hyperref}

\newlength{\bibitemsep}
\newlength{\bibparskip}
\let\oldthebibliography\thebibliography
\renewcommand\thebibliography[1]{%
  \oldthebibliography{#1}%
  \setlength{\parskip}{\bibitemsep}%
  \setlength{\itemsep}{\bibparskip}%
}

\begin{document}
\begin{center}
\Large \bf{A discrete approach to stochastic parametrization and dimensional reduction in nonlinear dynamics}
\end{center}

\begin{center}
Alexandre J. Chorin and Fei Lu\\
\vspace{3mm}
Department of Mathematics, University of California, Berkeley, CA 94720,~USA.\\
Lawrence Berkeley National Laboratory, Berkeley, CA 94720,~USA.
\end{center}

\vspace{1mm}

{\bf Abstract} 
Many physical systems are described by nonlinear differential equations that are too complicated to solve in full. A natural way to proceed is to divide the variables into those that are of direct interest and those that are not, formulate solvable approximate equations for the variables of greater interest, and use data and statistical methods to account for the impact of the other variables. In the present paper the problem is considered in a fully discrete-time setting, which simplifies both the analysis of the data and the numerical algorithms. The resulting time series are identified by a NARMAX (nonlinear autoregression moving average with exogenous input) representation familiar from engineering practice. The connections with the Mori-Zwanzig formalism of statistical physics are discussed, as well as an application to the Lorenz 96 system.

\section{Introduction and outline}

There are many problems in science where the equations of motion are too complex
 for full solution, either because the equations are not certain or because the computational cost is too high, but one is interested only in the dynamics of a subset of the variables. Such problems appear, for example, in weather and climate modeling, see e.g. \cite{Kal03,Wil11}, in economics, see e.g. \cite{ZW13}, in statistical mechanics, see e.g.\cite{EM90,CH13}, and in the mechanics of turbulent flow, see e.g. \cite{BW02,Cho94}. In these settings it is
natural to look for simpler equations that involve only the variables of interest, and then use data
to assess the effect of the remaining variables on the variables of interest in some approximate way. In the present paper we focus on stochastic methods for doing this,
with data coming either from experimental observations or
from prior numerical calculations.

Consider a set of differential equations of the form:
\begin{equation}
\frac {d}{dt}x=R(x,y),\,\,\,\, \frac{d}{dt}y=S(x,y),
\label{main}
\end{equation}
where $t$ is the time, $x=(x_1,x_2,\dots,x_m)$ is the vector of resolved variables, and $y=(y_1,y_2,\dots,y_{\ell})$ is the vector of unresolved variables, with initial data $x(0)=\alpha,\,y(0)=\beta.$ Consider a situation where this system is too complicated to solve, but where data are available, usually as
sequences of measured values of $x$, assumed here to be observed with negligible observation errors.
One can write $R(x,y)$ in the form
\begin{equation}
R(x,y)=R_0(x)+z(x,y),
\label{reduce}
\end{equation}
where $R_0$ approximates $R(x,y)$ in some sense and is such that one is able to solve the equation
\begin{equation}
\frac {d}{dt}x=R_0(x).
\label{Rzero}
\end{equation}
The remainder $z(x,y)=R(x,y)-R_0(x)$, often called the ``unresolved tendency", is
 the contribution of the unresolved variables to the equation for $x$. Recent work has shown that
 $z$ can represent model error \cite{Har13,BH14} or model noise \cite{CH14};
 we call $z$ simply the ``noise".

A usual approach to the problem of computing $x$ is to identify $z$ from from data \cite{PS07,Wil05,CVE08}, i.e., find a concise approximate representation $\hat{z}$ of $z$ that
can be evaluated on the computer, and then solve the equation
\begin{equation}
\label{usual}
\frac {d}{dt}x=R_0(x)+\hat{z}.
\end{equation}
When $\hat{z}$ is a stochastic process, this is a ``stochastic parametrization". Equation (\ref{usual}) is an instance of a dimensionally-reduced equation, in the sense that it has fewer variables than equations (\ref{main}). However, this approach has some difficulties. In general the available data are measurements of $x$, not of $z$;
to find $z$ so that it can be identified one has to use equation (\ref{reduce}), and in particular
differentiate $x$ numerically, which is generally impractical or inaccurate because $z$ may
have high-frequency components or fail to be sufficiently smooth, and the data may not be available at sufficiently small time intervals (an illuminating analysis in a special case can be found in \cite{PSW09,ST12}).  If one can successfully identify $z$, equation (\ref{usual}) generally becomes a nonlinear stochastic differential system, where in general $\hat{z}$ at a given time depends on earlier values of $x$ and $\hat{z}$ (see the next section), that may be hard to solve with sufficient accuracy (see e.g \cite{KP99,MT04}).

Here we take a different approach. We note that equations (\ref{Rzero}) and (\ref{usual}) are always solved on the computer, i.e., in discrete form, that the data are always given at a discrete
collection of points, and that one wishes to determine $x$ but in general one is not interested in determining $z$. We can therefore avoid the difficult detour through a continuous $z$ followed by a discretization, by working entirely in the discrete setting as follows. We can pick once and for all a particular accurate discretization of equation (\ref{Rzero}) with a particular time step $\delta$,
\begin{equation*}
x^{n+1}=x^n+\delta R_{\delta}(x^n),
\end{equation*}
where $R_{\delta}$ is obtained, for example, from by a Runge--Kutta method, and where $n$ indexes the result after $n$ steps. As in the continuous case, this equation has to be corrected to account for the impact of the continuous variables, and here also for the possible inaccuracy of the numerical scheme. We use the data to identify the discrepancy sequence, $z_{\delta}^{n+1}=(x^{n+1}-x^n)/\delta-R_{\delta}(x^n)$, which are available from $x$ data without approximation. This is equivalent to identifying the following reduced system
\begin{equation} 
\label{disc}
x^{n+1}=x^n+\delta R_{\delta}(x^n)+\delta z_{\delta}^{n+1}, 
\end{equation}
which constitutes a discrete stochastic parametrization.

We assume, as one does in the continuous case, that the system under consideration is ergodic, so that its long-time statistics are stationary. The sequence $z_{\delta}^n$ becomes a stationary time series, which we represent by the widely used NARMAX (nonlinear auto-regression moving average with exogenous inputs) representation, with $x$ as an exogenous input. This representation makes it possible to integrate the numerical scheme into the reduced equations, and to take into account efficiently the non-Markovian features of the reduced system as well as model and numerical errors. There is no stochastic differential system to solve. It is important to note that identifying $z_{\delta}$ can be very different from identifying the continuous $z$. The question, in what sense does $z_{\delta}$ approximate $z$, is not relevant, since the goal is to calculate $x$ and $z_{\delta}$ is sufficient for the purpose. Note that $z_{\delta}$ should be a good approximation of $z$ whenever equation (\ref{Rzero}) can be effectively approximated by the first-order Euler scheme. For practical purposes, the discrete stochastic parametrization accomplishes everything that would be accomplished by a successful continuous parametrization followed by an accurate approximation. From now on, we drop the hat that distinguishes between a process and its identification and the subscript $\delta$ in $z_{\delta}$.

The rest of the paper is organized as follows. In the next section we summarize the
NARMAX methodology we use, make comparisons with earlier work, and explain the connections with the Mori-Zwanzig formalism of statistical physics. In section 3 we present our test problem, the Lorenz 96 equations.
In section 4 we present numerical results. In the final section we draw conclusions.

\section{The NARMAX representation}

We represent $z$ in the reduced system (\ref{disc}) by a variant of the NARMAX representation, see e.g. \cite{Bil13, BD02, FY03, Ham94}. Its advantages are its generality, and the extensive experience in its use. Most of the earlier
identifications can be viewed as special case of this representation, as
explained below. The greater generality of the NARMAX approach will be
increasingly important as model reduction methods are applied to more
complex problems.
Equation (\ref{disc}) becomes
\begin{equation}
\label{full}
x^{n+1}=x^n+\delta R_{\delta}(x^n)+\delta z^{n+1},\,\,\,z^{n+1}=\Phi^{n+1}+\xi^{n+1},
\end{equation}
for $n=1,2,\dots$,
where the $\xi^{n+1}$ are independent Gaussian random variables with mean zero and variance $\sigma^2$, and
$\Phi^n$ is the sum:
\begin{equation}
\Phi ^{n}=\mu
+\sum_{j=1}^{p}a_{j}z^{n-j}+
\sum_{j=1}^{r}\sum_{i=1}^{s}b_{i,j}P_{i}(x^{n-j})+\sum_{j=1}^{q}c_{j}\xi^{n-j},
\label{NARMAX}
\end{equation}
$\mu $, $\sigma ^{2}$ and $\left\{
a_{j},b_{i,j},c_{j}\right\}$ are parameters to be inferred
from data, and the exogenous inputs $P_{i},\ i=1,\dots,s$ are functions to be determined;  to simplify the notations, equation (\ref{NARMAX}) has been written as if equation (\ref{full}) were scalar. This is the NARMAX representation of $x$ and $z$. In equation (\ref{NARMAX}), the terms in $z$ are the autoregression part of order $p$, the terms in $\xi$ are the moving average part of order $q$, and the $P_{i}(x)$ are the exogenous inputs.
Note that in the reduced system (\ref{full}), the $z^n$ term can be eliminated; indeed, if we substitute
$z^n=(x^n-x^{n-1})/\delta-R_{\delta}(x^{n-1})$ into (\ref{NARMAX}) and the second of equation in (\ref{full}), we obtain a NARMA representation for $x$. A suitable choice of the functions $P_i$ will be discussed in the context of a specific example and will connect the representation to the approximation of equation (\ref{Rzero}).

To implement the NARMAX representation, one has to determine its structure and estimate the parameters. While NARMAX has been widely studied (see e.g. \cite{Bil13, DC05, Chen10, Han76,Sto86} and references therein), one should use the earlier work with caution, especially in the detection of structure by least-squares-based methods, because in the standard NARMAX, unlike here, the exogenous process is independent of the noise process. Suppose one has selected a structure, that is, chosen the functions $P_{i}$ and the orders $\left( p,r,s,q\right) $ in the representation (\ref{NARMAX}). Since the representation is linear in the parameters $\theta =\left(\mu, a_{j},b_{i,j},c_{j}\right) $, these parameters can be estimated using conditional likelihoods as follows. The sequence $%
\left\{ z^{n}\right\}$ for $n=1,2,\dots,N$ can be computed from the data using the first equation in (%
\ref{full}). Once the values of $\xi^{1},\dots ,\xi^{q}$ are
known, the noise sequence $\xi^{n}$ for $q+1 \le n \le N$ can be
computed from $\xi^{n}=z^n -\Phi^{n}$. This leads to the conditional log--likelihood of the observations $x^n$ for
$q+1 \le n \le N$ (up to a constant):
\begin{equation*}
l(\mathbf{\theta }|\xi^{1},\dots ,\xi^{q})=-\sum_{n=q+1}^{N}\frac{%
\left\vert z^n-\Phi ^{n}\right\vert ^{2}}{2\sigma ^{2}}-\frac{%
N-q}{2}\ln \sigma ^{2}.
\end{equation*}%
Since the NARMAX representation is linear in its parameters other than $\sigma^2$, the
log-likelihood $l(\theta |\xi^{1},\dots ,\xi^{q})$ is quadratic in these
parameters, its gradient can be easily computed and the maximum
likelihood estimator (MLE) can be obtained by standard gradient-based
optimization methods, such as the quasi-Newton method. If the reduced system is ergodic,
the MLE is asymptotically consistent, and the initial values of $\xi^{1},\dots ,\xi^{q}$ do not affect the result (see e.g. \cite{Ham94, Han76}). For convenience, we set $\xi^{1}=\dots =\xi^{q}=E[\xi^{1}]=0$.

We remark that for the above Gaussian likelihood, the MLE is the same as
the least squares estimator for the parameters, which has been widely used
in control (see e.g. \cite{DC05, Chen10}). As shown in \cite{DC05}, the above estimation
procedure can be made recursive; also, the noise sequence need not be Gaussian.

The detection of the representation's structure, however, is less straightforward, as is well-known, see \cite{Bil13,FY03}.
Because in our problem the exogenous processes are not independent of the noise, popular techniques such as orthogonal least squares and error reduction ratios (see e.g. \cite{Bil13} and references therein), do not work.
We use the following criteria for selecting the structure of the representation: (1) it should fit the long-term statistics of the resolved variables, such as the mean and autocorrelation function; (2) as the size
of the data increase, the estimated parameters should converge; (3) the estimated
parameters should reflect features of the resolved variables, such as symmetry
properties. We find structures that satisfy these criteria by trial and error.

It is of interest to relate the NARMAX representation to the Mori-Zwanzig (MZ) formalism \cite{Zwan73,Zwan01,EM90,CH13}. The overall setting in the MZ representation is the same as here: one has data $\alpha$ for the $x$ variables, and one samples data $\beta$ for $y$ from a given initial pdf.
The MZ formalism creates the approximation $R_0(x)$ in equation (\ref{Rzero}) by conditional
expectation: $$R_0(x)=E[R(x,y)|x],$$ where $E[a|b]$ is the expected value of $a$ with respect to the initial
measure given $b$. When the system is ergodic and the initial pdf for $\beta$ is the invariant measure conditioned by $\alpha$, this is the best least-squares approximation of $R(x,y)$ by a function of $x$.
The MZ formalism then yields an expression for $z(x,y)$ in equation
(\ref{reduce}) as a sum of a noise and a non-Markovian memory/dissipation term, corresponding to $\xi^{n+1}$ and $\Phi^{n+1}$ in equation (\ref{NARMAX}); note that in (\ref{NARMAX}) $R_0$ is not restricted to the MZ recipe.
The MZ expressions are exact, and prove the need for the representation of $z$ to take the memory into
account by including information from earlier steps.

Once the initial data $y(0)=\beta$ have been sampled, the MZ equations are deterministic; the MZ formalism proposes to follow in time one particular sample path of the system. For a chaotic system such as the one discussed in the next section, this may not be computationally feasible. The representation here looks for sample paths of a stationary stochastic process whose statistics agree with the statistics defined
by the equations of motion. This is a related but less ambitious and more feasible task.

The evaluation of the various terms in the MZ formalism is difficult; as far as we know, there is
only one case in the literature \cite{CHK02} where it has been successfully carried
out in full for a nonlinear problem that is not completely trivial. The MZ formalism is a
useful starting point for analytic
approximations (see e.g. \cite{CHK02,Sti13})
but it is difficult to use it to construct reduced
models from data when suitable
analytic approximations are not available.
The formalism here is a more tractable way to use data for dimensional
reduction, and generalizes the MZ formalism
to a broader class of approximations.

There is a large literature on data-based dimensional reduction. In \cite{Wil05}, see also \cite{AMP13}, the noise $z$ is represented as the sum of an approximating polynomial in $x$ obtained by regression and a one-step autoregression; the details are in the next section where we compare our results to those in \cite{Wil05}. The shortcomings of this representation as a general tool are that it does not allow $z$ to remember past values of $x$, and that the autoregression term is not necessarily small, making it difficult to solve the equations accurately. Furthermore, in \cite{Wil05}, numerical values of a continuous $z$ are obtained by finite differences.
In \cite{CVE08,Kwa12} the noise is represented by a conditional Markov chain that depends on both current
and past values of $x$; the Markov chain is deduced from data by binning and counting, assuming that exact observations of $z$ are available. It should be noted that the Markov chain representation is intrinsically discrete, making this work close to ours in spirit. In \cite{MH13} the noise is treated as continuous and represented in a form that
is partly analogous to the NARMAX representation once one translates from the continuum to the grid.
An earlier construction of a reduced approximation can be found in \cite{CH14}, where the approach was not yet fully discrete. Other interesting related work can be found in
\cite{CKG11,KCG15, DN07,DD09}.

\section{The Lorenz 96 equations }

The Lorenz 96 equations \cite{Lor95} are a set of chaotic differential equations
that is often used as a metaphor for the atmosphere. It has been
widely used as a test bench for various dimensional reduction and stochastic parametrization methods
\cite{FVE04,Wil05, CVE08, Kwa12, AMP13}. Following \cite{CVE08},
we use the following formulation of the equations:
\begin{align*}
\frac{d}{dt}x_{k}&=x_{k-1}\left( x_{k+1}-x_{k-2}\right) -x_{k}+F+z_{k},
\label{L96_x} \\
\frac{d}{dt}y_{j,k}&=\frac{1}{\varepsilon }%
[y_{j+1,k}(y_{j-1,k}-y_{j+2,k})-y_{j,k}+h_{y}x_{k}]
\end{align*}
with
$z_{k}=\frac{h_{x}}{J}\sum_{j}y_{j,k}$,
and $k=1,\dots ,K$, $j=1,\dots ,J$.
The indices are cyclic,
$x_{k}=x_{k+K},~y_{j,k}=y_{j,k+K}$ and $y_{j+J,k}=y_{j,k+1}$. The
system is invariant under spatial translations, and the statistical
properties are identical for all $x_{k}$. The formulation here is
equivalent to the original formulation by Lorenz (see e.g. \cite%
{FVE04,Kwa12}); the parameter $\varepsilon$ measures the time-scale separation
between the resolved variables $x_{k}$ and the unresolved variables $y_{j,k}$. We set $\varepsilon
=0.5$, so that there is no significant time scale separation between resolved and
unresolved processes, as is both more realistic and more difficult to handle for
parameterizations (see \cite{FVE04} and references therein).
The other parameters are $K=18,J=20,F=10,h_{x}=-1$ and $h_{y}=1$. The ergodicity of the Lorenz 96 system has been numerically verified in earlier work (see e.g. \cite{FVE04}) and we do not revisit this issue.

In the following section we present numerical results produced by our NARMAX scheme and compare them to those in \cite{Wil05} labeled POLYAR.
We do not compare with the results  of \cite{CVE08} because they require exact observations of $z$.
In  \cite{Wil05}, the continuous $z$ is estimated by finite differences:
\begin{equation*}
z_{k}(t)\approx \frac{x_{k}(t+\delta )-x_{k}(t)}{\delta }-x_{k-1}\left(
x_{k+1}-x_{k-2}\right) +x_{k}-F.
\end{equation*}%
Then a polynomial regression of the form $
z_{k}(t)=P(x_{k}(t))+\eta _{k}(t)
$
is used to fit the data $\{x_{k}(n\delta),z_{k}(n\delta) \}$,
where $P(x)$ is an approximating polynomial obtained by least squares, and $\eta _{k}(t)$
is a one-step autoregression (AR(1)) with parameters estimated from the time series $
z_{k}(n\delta )-P(x_{k}(n\delta)),\,\,$ for $n=1,2,\dots.$. This
leads to the following reduced stochastic equation:
\begin{equation}
\frac{d}{dt}x_{k}=x_{k-1}\left( x_{k+1}-x_{k-2}\right)
-x_{k}+F+P(x_{k})+\eta _{k},
\label{x_POLYAR}
\end{equation}
where $\eta _{k}$ is an autoregression of the form
\begin{equation}
\eta _{k}(t+\delta)=\phi\eta _{k}(t)+\sigma \xi _{k}(t)  \label{AR}
\end{equation}
where $\phi,\,\sigma$ are constants deduced from the data, and the
 $\xi _{k}(t)$ are independent identically distributed
Gaussian random variable with mean zero and variance one, for each component $k=1,\dots ,K$ of the
equation.
This reduced system is solved as follows: given the current time vectors $({\eta}_k(t)\,,x_k(t))$, the next time-step ${\eta}_k(t+\delta )$ is calculated from (\ref{AR}),
and then ${x_k}(t+\delta )$ is computed by integrating (\ref{x_POLYAR})
by a fourth-order Runge--Kutta methods, with ${\eta}(t)$ kept constant during each time step.

In the NARMAX scheme, we use the representation (\ref{NARMAX}), choosing one of the functions $P_i(x)$ to be $R_{\delta}(x)$ from the approximation (\ref{Rzero}) and the others to be powers of x. This connects the numerical scheme with the representation of the noise. We select the structure and estimate
the parameters as described earlier.
The parameters are the same for each spatial component, refecting the spatial symmetry of the equation. Each component of $z$ remembers only its
own past and the past of the corresponing component of $x$. The term $\Phi^n$ in equation (\ref{NARMAX})becomes:
\begin{align}
\Phi^n& =\mu
+\sum_{j=1}^{p}a_{j}z^{n-j}
+\sum_{j=1}^{r}\sum_{l=1}^{d_x}b_{j,l}(x^{n-j})^{l}  \notag \\  
& +\sum_{j=1}^{s}\sum_{l=1}^{d_R}c_{j,l} ( R_{\delta}(x^{n-j}))^{l}
+\sum_{j=1}^{q} d_{j} \xi^{n-j}.
\label{NARMAX_L96}
\end{align}
The determination of the numerical parameters in this representation is part of the
calculation and the values used are listed in the next section.

\section{Numerical results}
In the numerical runs, we generate data for $x_{k}$ by integrating the full Lorenz 96 system with parameters $\left(\varepsilon ,K,J,F,h_{x},h_{y}\right) =(0.5,18,20,10,-1,1)$, using a fourth order Runge-Kutta method with time step $0.001$. We consider two cases: one in which the observations are made at intervals of $\delta=0.01$
and one at which they are made at intervals $\delta=0.05$; the first case corresponds to a case discussed in (\cite{CVE08,Kwa12}); in the second case, the data are slightly sparser. In each case, we make observations at
at $N=5\times 10^5$ points; this requires that the full system be integrated for $5000$ and $25000$ time units,
respectively.

\begin{table}[tbp]
\caption{The estimated parameters in the POLYAR system; the column labeled $j$ for $j=0,\dots,5$ contains the coefficient of $x$ to the power $j$.}
\label{tab1}\centering
\begin{tabular}{ccccccccc}
\hline
& $5$ & $4$ & $3$ & $2$ & $1$ & $0$ & $\phi$ & $\sigma$
 \\ \hline
$\delta =0.01$ & $-0.00002$ & $0.0004$ & $-0.0002$ & $-0.0258$ & $-0.3567$ & $0.0529$ & $0.9948$ & $0.9397$ \\
$\delta =0.05$ & $-0.00003$ & $0.0009$ & $-0.0035$ & $-0.0137$ & $-1.0030$ & $1.3919$ & $0.7626$ & $11.3857$ \\ \hline
\end{tabular} %
\end{table}

\begin{table}[tbp]
\caption{The estimated parameters in the NARMAX model.}
\label{tab2}\centering
\begin{tabular}{cccccccc}
\hline
& $a_{1}$ & $b_{1,1}$ & $b_{1,2}$ & $d_{1}$ & $\, $ & $ \mu $ & {$\sigma ^{2}$} \\ \hline
$\delta =0.01$ & $0.9782$ & $-0.1271$ & $0.1132$ & $0.9997$ & \, & $0.0115$ & ${0.0004}$ \\
& $a_{1}$ & $b_{1,1}$ & $b_{1,2}$ & $b_{1,3}$ & $c_{1,1}$ & $\mu $ & $\sigma^{2}$ \\ \hline
$\delta =0.05$ & $0.8879$ & $-0.0712$ & $-0.0002$ & $0.0002$ & $-0.0084$ & $0.0556$ & $0.0284$ \\ \hline
\end{tabular}
\end{table}%

\begin{figure}
\centering
\begin{subfigure}[b]{0.4\textwidth} 
                \includegraphics[width=\textwidth]{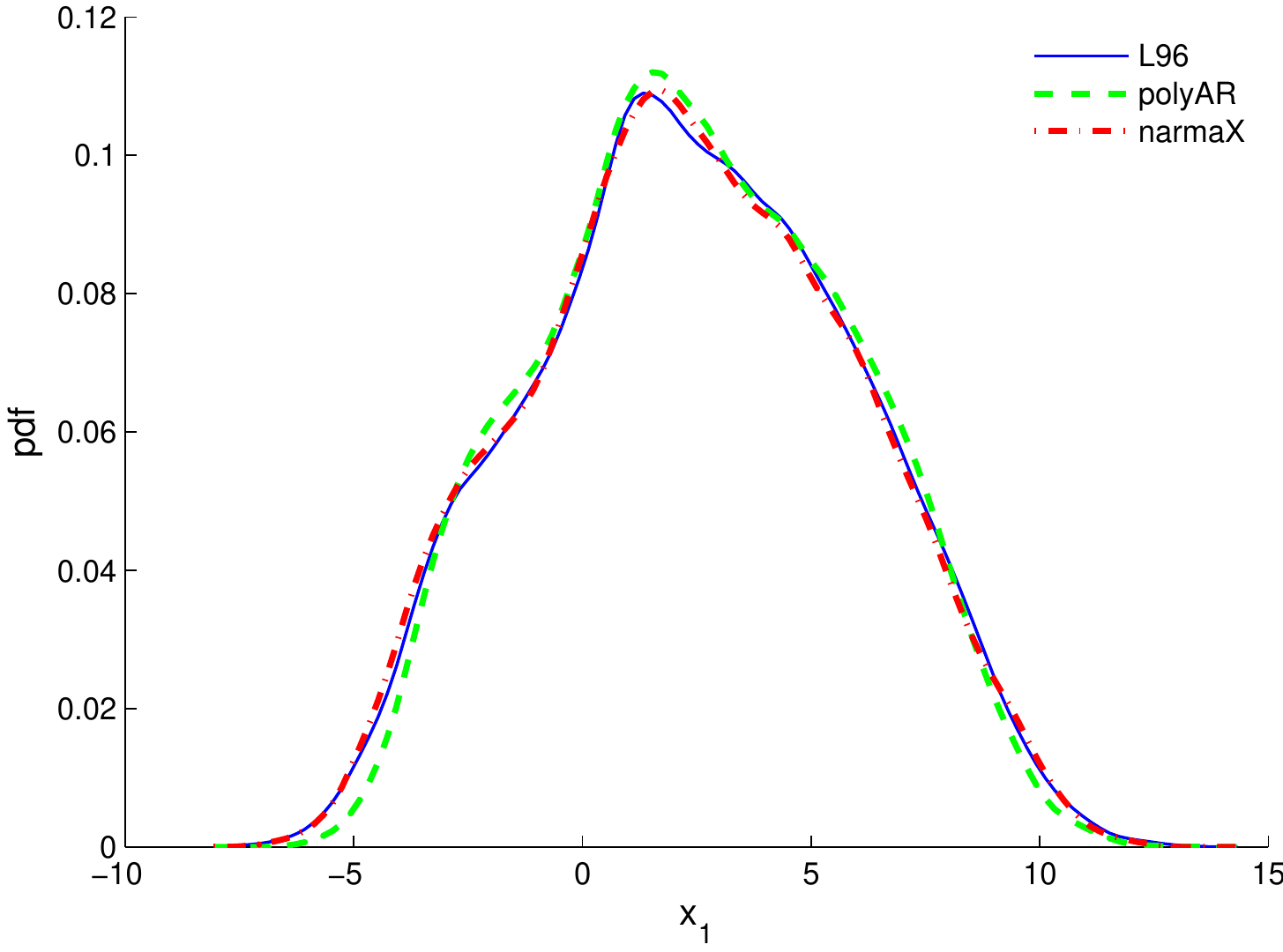}
        \end{subfigure}
\begin{subfigure}[b]{0.4\textwidth}
                \includegraphics[width=\textwidth]{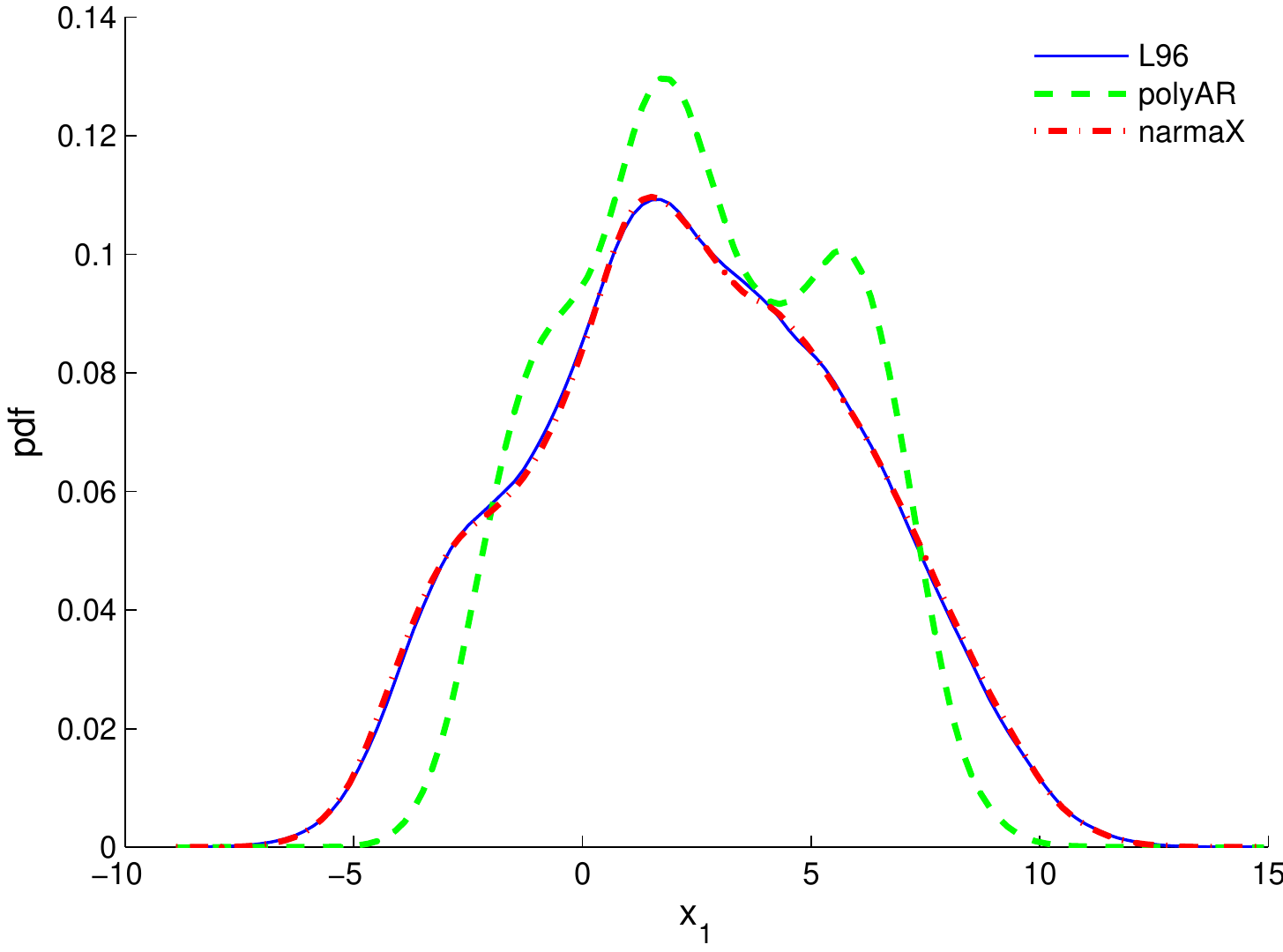}
        \end{subfigure}
\begin{subfigure}[b]{0.4\textwidth}
                \includegraphics[width=\textwidth]{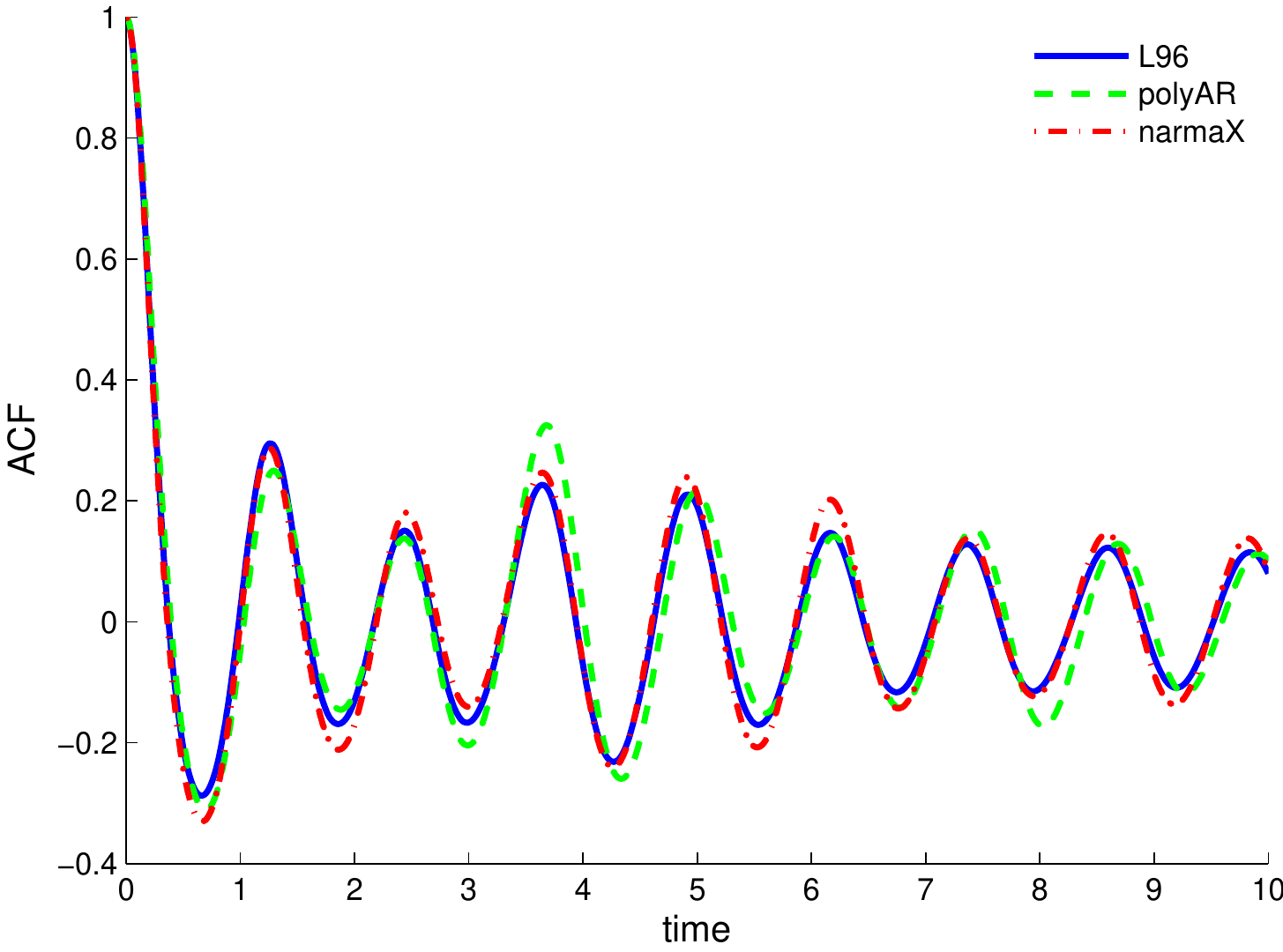}
          \end{subfigure}
\begin{subfigure}[b]{0.4\textwidth}
                \includegraphics[width=\textwidth]{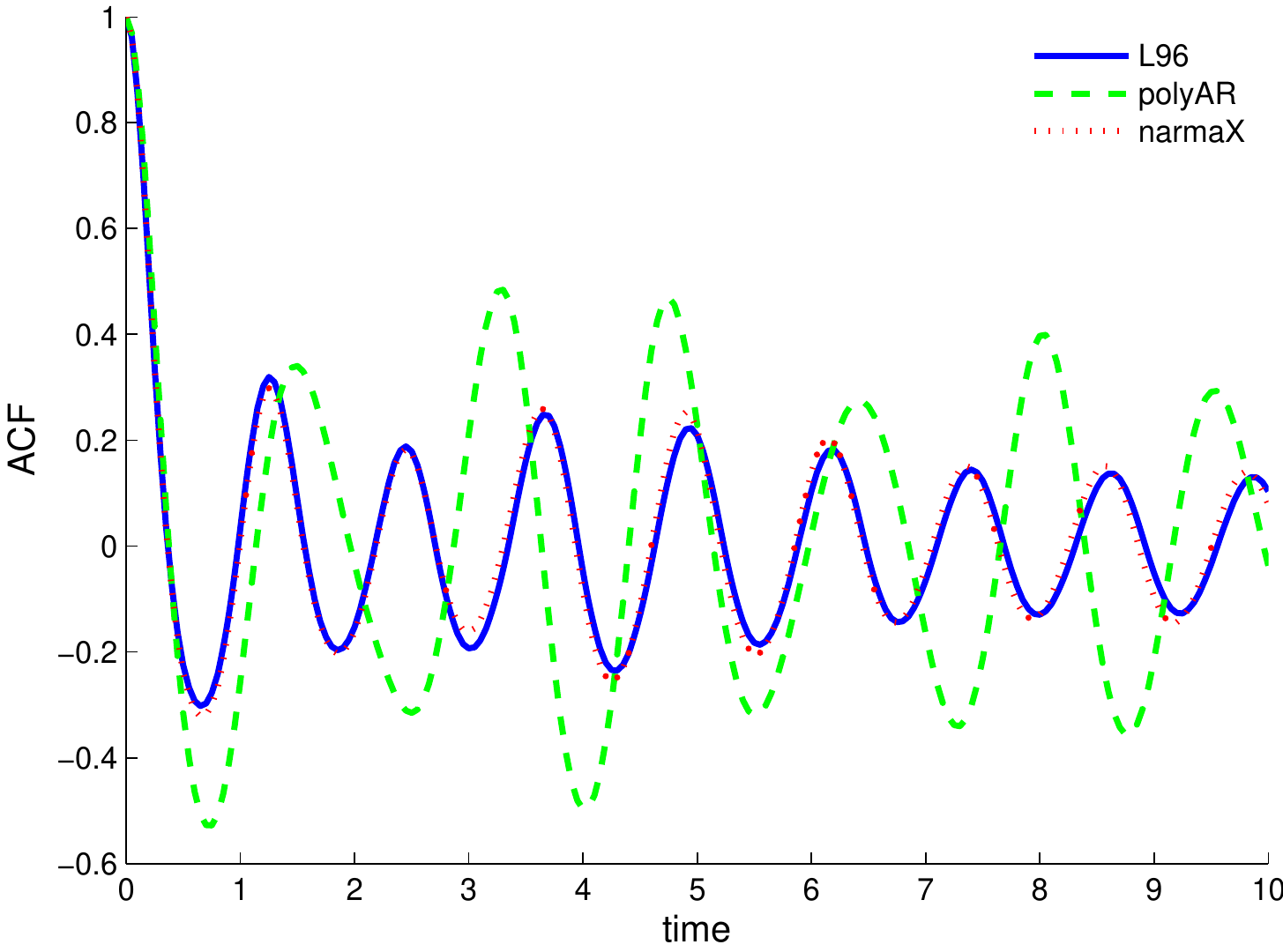}
         \end{subfigure}
\begin{subfigure}[b]{0.4\textwidth}
                \includegraphics[width=\textwidth]{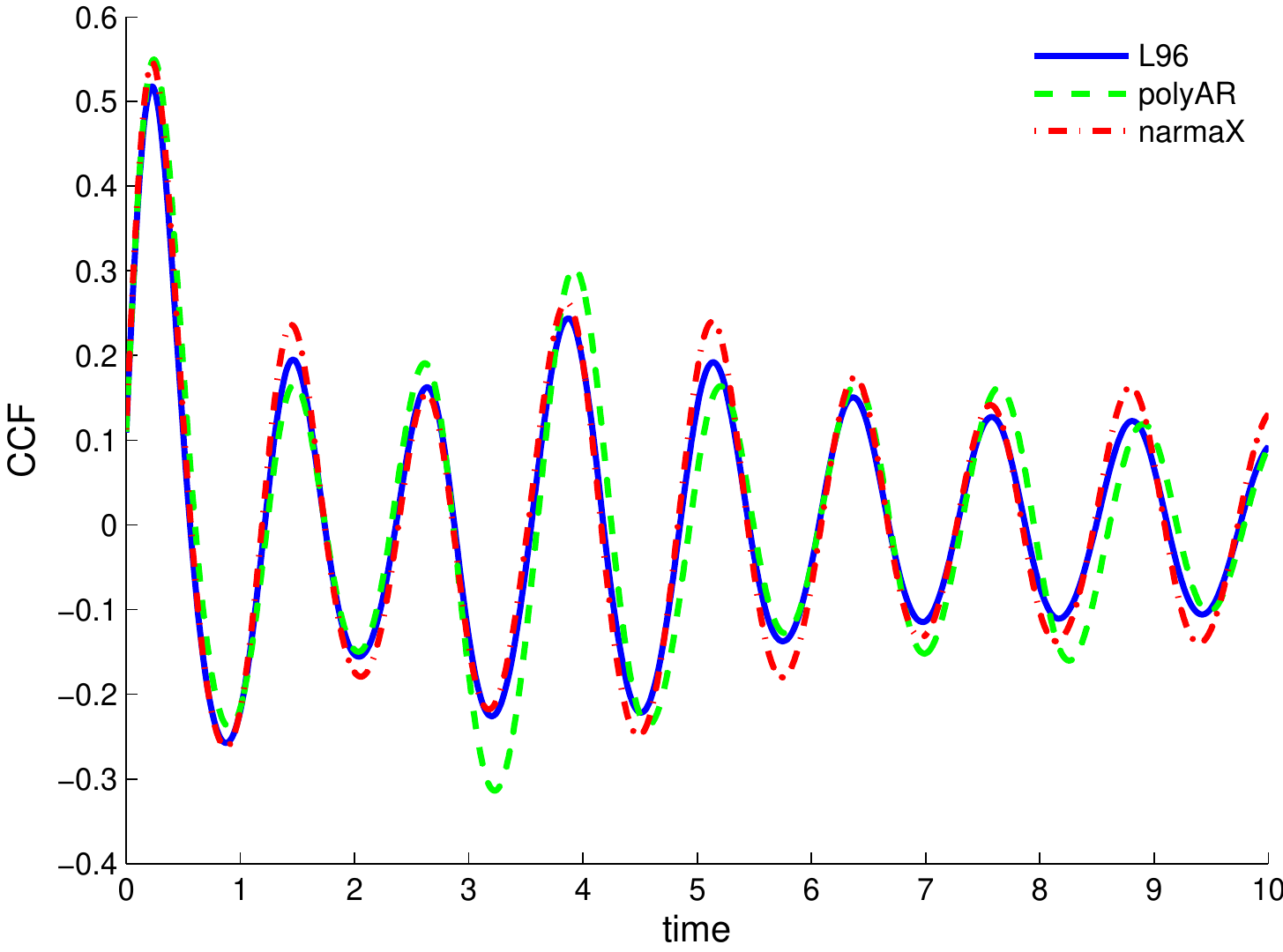}
         \end{subfigure}
\begin{subfigure}[b]{0.4\textwidth}
                \includegraphics[width=\textwidth]{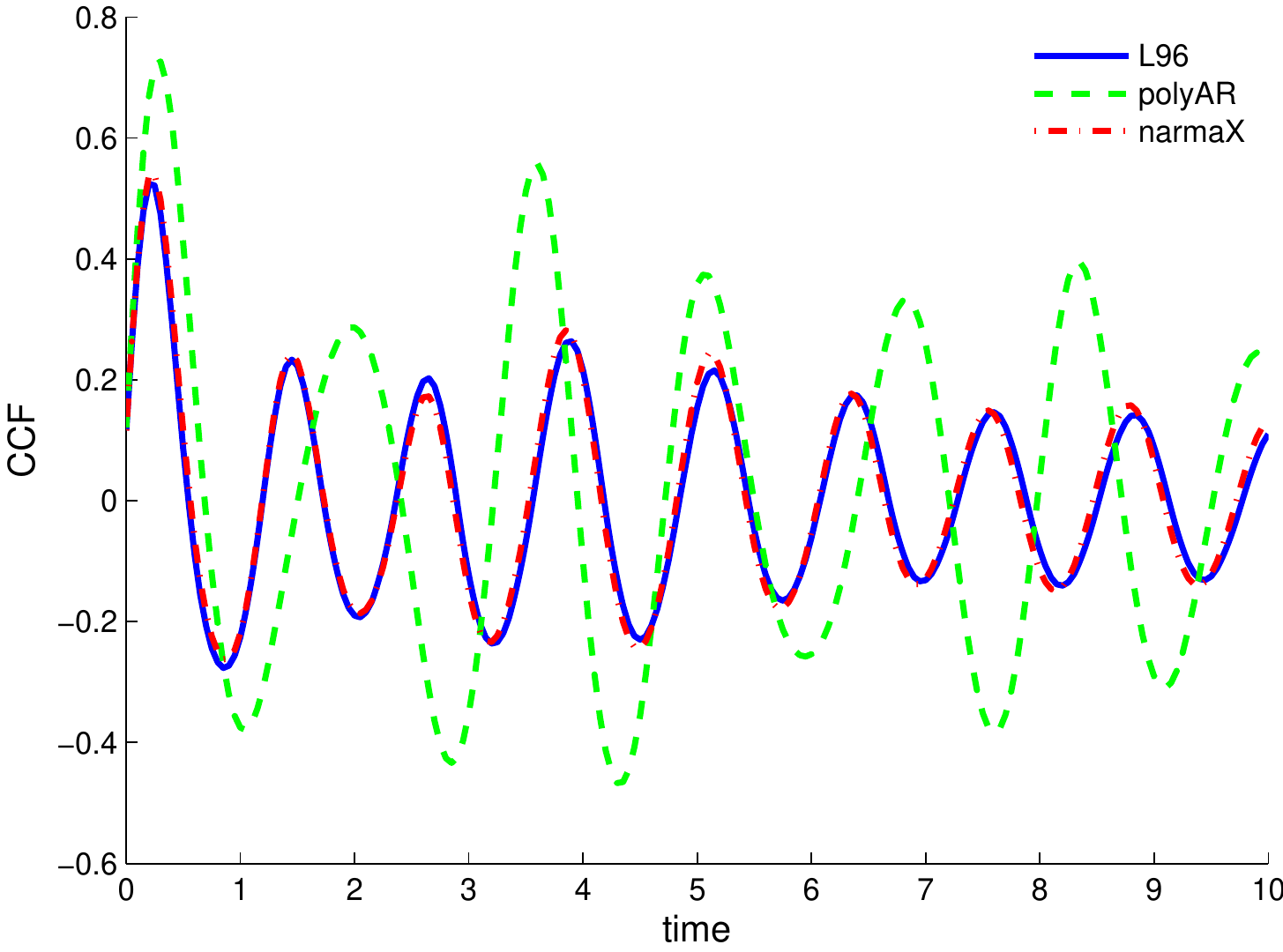}
        \end{subfigure}
\caption{Probability density functions (pdf), autocorrelation functions (ACF) and cross-correlation functions (CCF) of $x_1$, produced by the full Lorenz 96 model, the POLYAR system and the NARMAX system. The left column is the case $\delta=0.01$, and the right column is the case $\delta=0.05$.}
\label{fig:pdfACF}
\end{figure}

For the POLYAR equations of \cite{Wil05}, a fifth-order polynomial is used for both observation settings. Increasing the order further produces small coefficients for the higher degree terms, which do not reduce the variance of noise. The estimated parameters, i.e. the coefficients of the polynomial and the parameters for the autoregression, for the first component $x_1$ of $x$ are shown in Table 1.

For the NARMAX equation (\ref{NARMAX_L96}), we estimated $(p,r,s,q)=(1,2,0,1)$ and $\left(
d_{x},d_{R}\right) =\left( 1,0\right) $ for the case $%
\delta =0.01$, and for the case $\delta =0.05$,  $(p,r,s,q)=(1,1,1,0)$
and $\left( d_{x},d_{R}\right) =\left( 3,1\right) $. The
estimated parameters for $x_{1}$ are shown in Table 2.

First, we compare the statistics of the solutions generated by the two reduced systems with the statistics of the full system. We integrate the reduced equations in both cases and obtain values at $5\times 10^5 $ points. We calculate the following quantities from the reduced equations as
well as from the full Lorenz 96 equations: the mean, the standard deviation, the
probability density function (pdf), the Kolmogorov--Smirnov statistic,
the autocorrelation function (ACF) of $x_1$, and the cross-correlation function (CCF) between $x_1$ and $x_{2}$.
The pdf of $x_1$ for the full Lorenz 96 system is well reproduced by
both reduced systems when $\delta=0.01$, see left top in Figure~\ref{fig:pdfACF}.
In the sparser data case $\delta=0.05$ the NARMAX
 equations produce a much better pdf than the POLYAR equations, see right top in Figure~\ref{fig:pdfACF}.
Table 3 displays the mean, the standard deviation,
and the Kolmogorov--Smirnov statistic that compare the cumulative
distributions of the full Lorenz 96 system and with that of the reduced equations. The NARMAX system
is more accurate than the POLYAR system in both cases.
The autocorrelation function and the cross correlation function are well reproduced by both reduced systems when $\delta =0.01$, see left middle and left bottom of Figure~\ref{fig:pdfACF}. When $\delta =0.05$, however, the autocorrelations and the cross correlations reproduced by the POLYAR miss the
amplitude of oscillation and exhibit a phase shift from those of the full Lorenz 96 equations, while the NARMAX system remains accurate, see right middle and right bottom of Figure \ref{fig:pdfACF}.

\begin{table}[tbp]
\caption{The mean, standard deviation and the Kolmogorov--Smirnov statistic (D).} 
\label{tab3}\centering
\begin{tabular}{cccccccc}
\hline
& \multicolumn{3}{c}{$\delta =0.01$} &  & \multicolumn{3}{c}{$\delta =0.05$}
\\ \cline{2-4}\cline{6-8}
& mean & std. & D &  & mean & std. & D \\ \hline
full Lorenz 96 & $2.4506$ & $3.5230$ &  &  & $2.3978$ & $3.5222$ &  \\
POLYAR & $2.5335$ & $3.3807$ & $0.0183$ &  & $2.6031$ & $2.8564$ & $0.0747$ \\
NARMAX & $2.4113$ & $3.5270$ & $0.0055$ &  & $2.4293$ & $3.5402$ & $0.0049$ \\ \hline
\end{tabular}%
\end{table}

We now investigate how well a reduced system predicts the behavior of the full system
by calculating mean trajectories of the reduced systems and comparing them
with a true trajectory of the full Lorenz 96 system, as follows.
First we integrate the full Lorenz 96 system for $10\times N_0$ time units, and store the results as $N_0$ short trajectories of 10 time units each. For each short true trajectory, we perform $N_{ens}$ integrations of the reduced systems over 10 time units, starting all ensemble members from the same several-step initial conditions as the corresponding full solution -- several initial steps are needed to initialize $\eta$ in POLYAR and $\xi$ in NARMAX. We do not introduce artificial perturbations into the initial conditions, because the exact initial conditions for $x$ are known, and by initializing $\eta$ and $\xi$ from data, we preserve the memory of the system so as to generate better ensemble trajectories. We then average the solutions of the reduced equations in each subinterval and compare these averages with the trajectories of the full system
by calculating the root-mean-square error (RMSE) and anomaly correlation (ANCR) between them;
the former measures the average difference between trajectories while
the latter measures the average correlation between them; the formulas and their rationale can be found e.g. in \cite{CVE08}.

\begin{figure}[tbp]
\centering
\begin{subfigure}[b]{0.4\textwidth}
                \includegraphics[width=\textwidth]{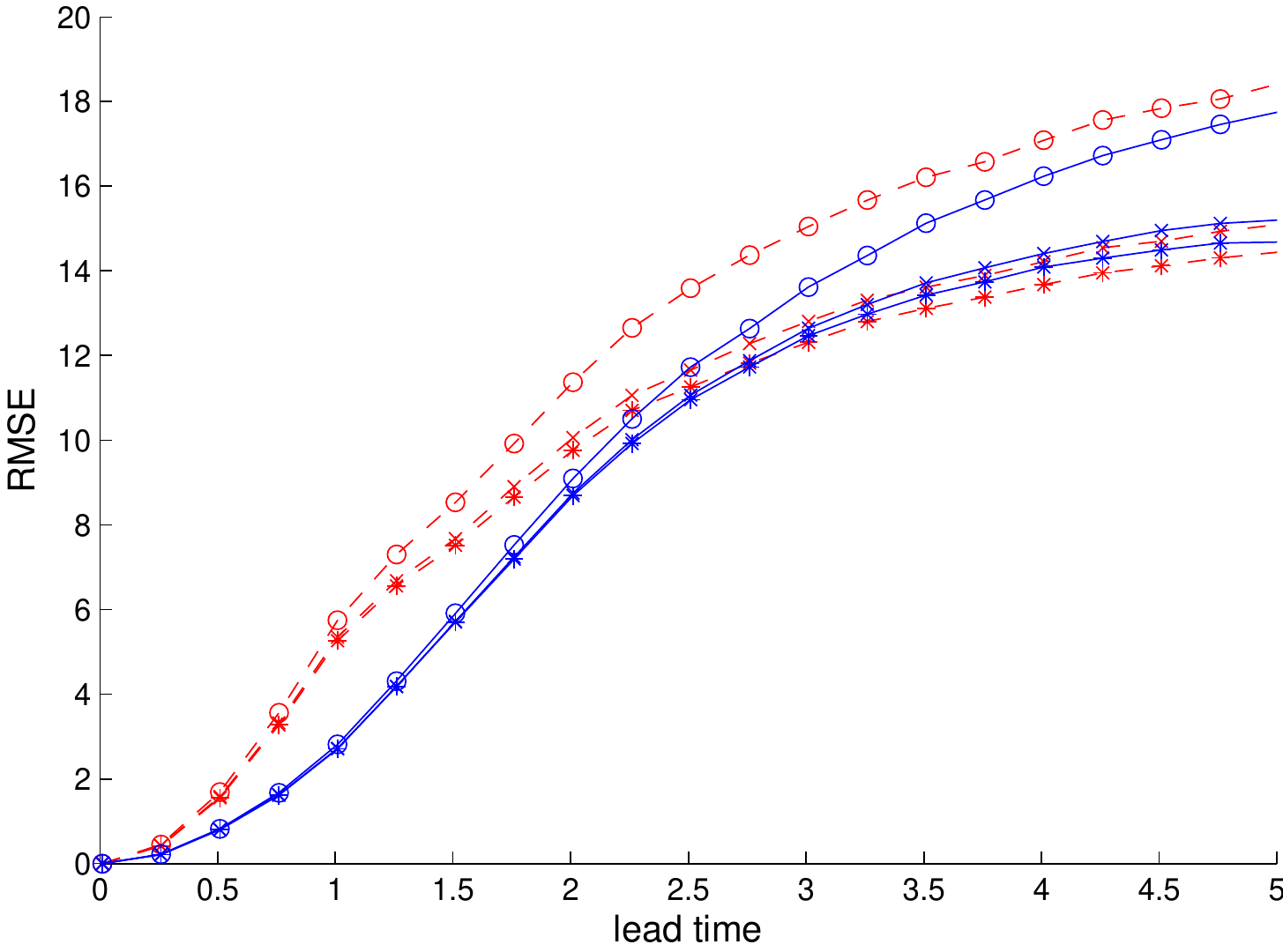}
        \end{subfigure}
\begin{subfigure}[b]{0.4\textwidth}
                \includegraphics[width=\textwidth]{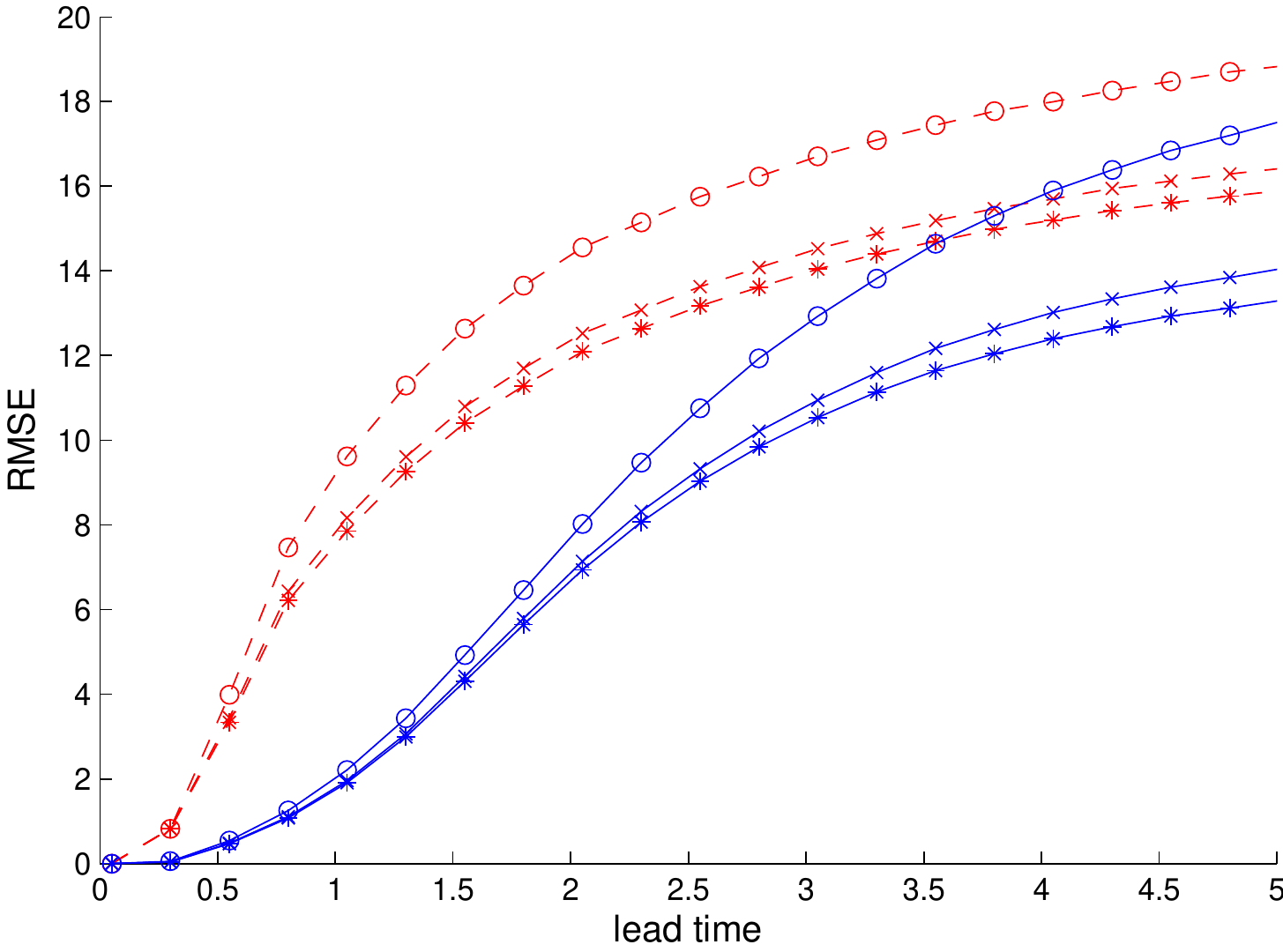}
        \end{subfigure}
\begin{subfigure}[b]{0.4\textwidth}
                \includegraphics[width=\textwidth]{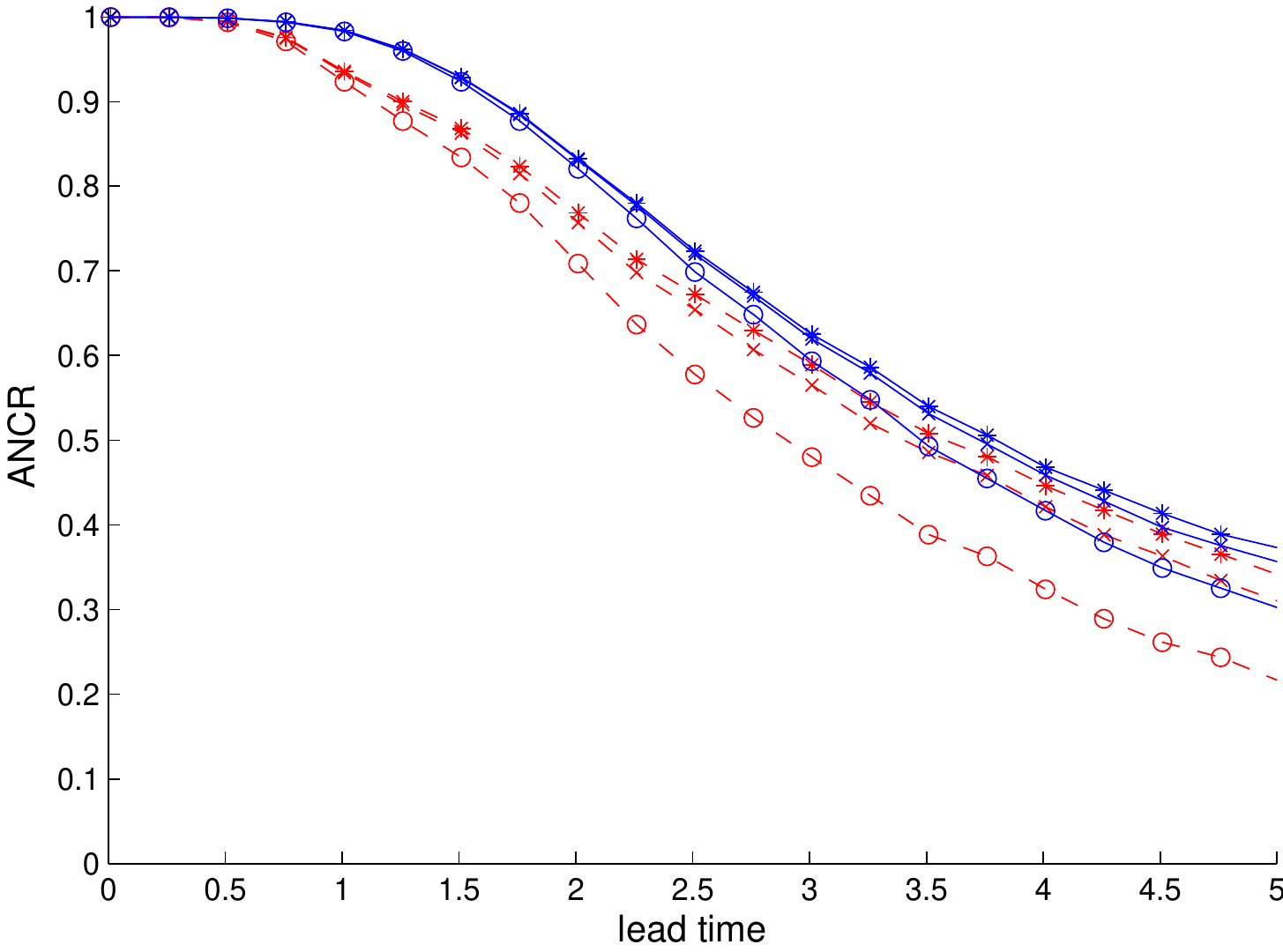}
        \end{subfigure}
\begin{subfigure}[b]{0.4\textwidth}
                \includegraphics[width=\textwidth]{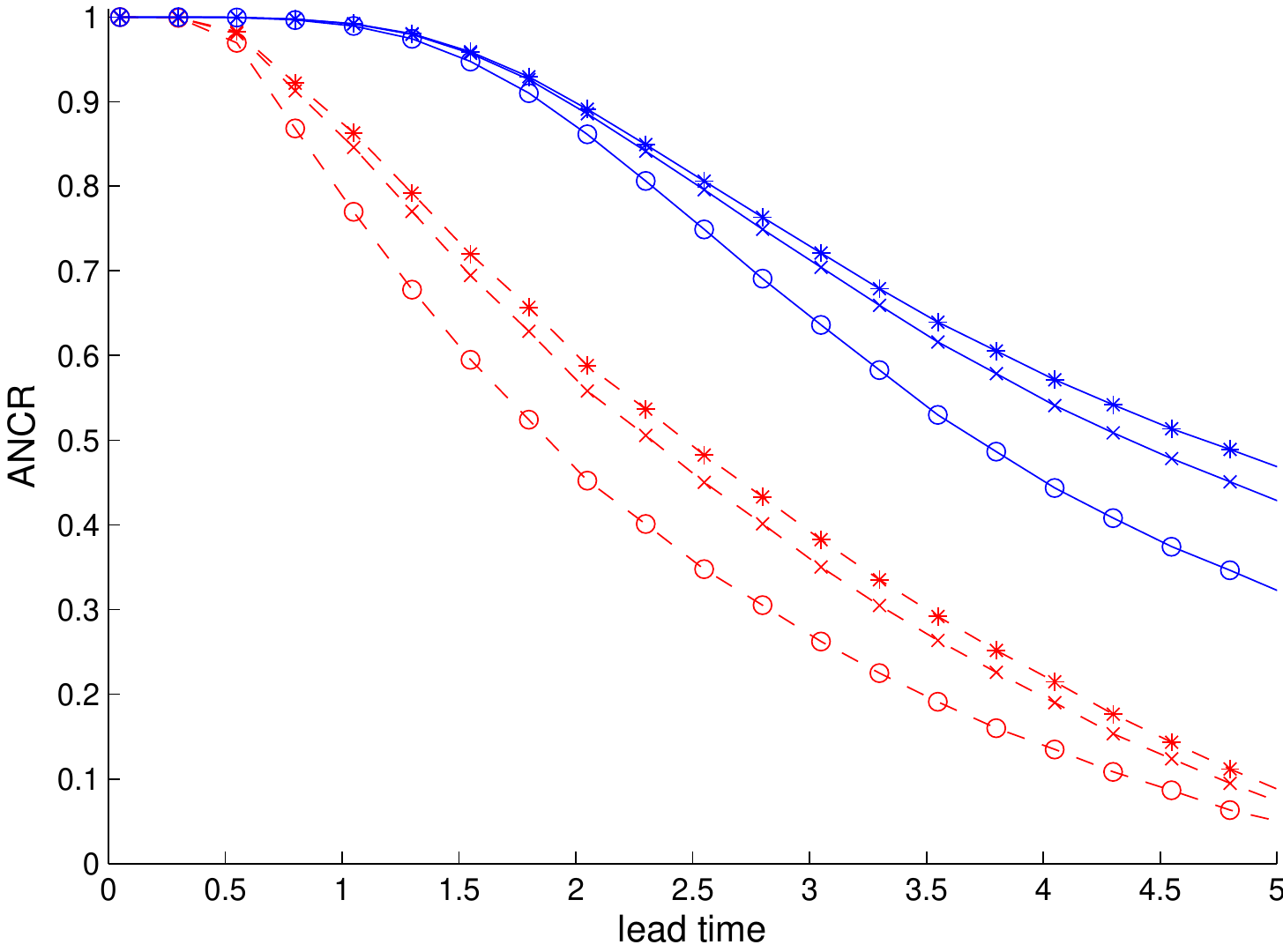}
        \end{subfigure}
\caption{Root mean square error (RMSE) and anomaly correlation (ANCR) of ensemble forecasting, produced by the NARMAX system (solid line) and the POLYAR system (dashed line), for different ensemble sizes: $N_{ens}=1$ (circle marker), $N_{ens}=5$ (cross marker), and $N_{ens}=20$ (asterisk marker). The left column is the case $\delta=0.01$, and the right column is the case $\delta=0.05$.}
\label{fig:RMSEancr}
\end{figure}

Results for RMSE and ANCR are shown in Figure~\ref{fig:RMSEancr}, where we use $N_0=10000,\, N_{ens}=1,5,20$ and the number of steps where initial conditions are imposed is $n_0=\max \{ 1,p,r,s,2q\}+1$. In the case $\delta=0.01$, the NARMAX reduction performs slightly better than the POLYAR reduction. In the case $\delta=0.05$, the NARMAX reduction provides a significant improvement over the POLYAR reduction. For example, the forecast lead time at which the anomaly correlation drops below 0.6 is extended from $\tau=2$ to $\tau=4$ in the case $N_{ens}=20$.

\section{Conclusions}
We have presented a discrete approach to data-based dimensional reduction and stochastic parametrization, in which the problem is consistently treated as discrete, obviating earlier difficulties in estimating noise from measurements and in approximating reduced continuum equations. Within this discrete approach, we have identified the reduced system via the NARMAX representation. This generalizes earlier work, in particular by making it easy to include memory effects in full. We have tested the resulting algorithm on the Lorenz 96 system of equations, often used as a test bench for dimensional reduction schemes; our construction did better than earlier work in reproducing the dynamics and the long-range statistics of the variables of interest, most conspicuously in a problem where the data were sparse. We related our representation to the Mori-Zwanzig formalism and suggested that our methods can be used to construct data-based implementations of this formalism. We expect the advantages of our modeling to become even more marked as it is applied to increasingly complex problems.

\section{Acknowledgements}
This work was supported in part by the Director, Office of Science, Computational and Technology Research, U.S. Department of Energy, under Contract No. DE-AC02-05CH11231, and by the National Science
Foundation under grants DMS-1217065 and DMS-1419044. The authors would like to thank their Berkeley colleague Dr. Matthias Morzfeld, Prof. Kevin Lin of the University of Arizona, Prof. Xuemin Tu of the University of Kansas, and Prof. Robert Miller of the State University of Oregon for helpful comments and good advice.
\bibliographystyle{unsrt}



\end{document}